\documentclass{fundam}
\usepackage{eufrak}
\usepackage{euler}
\usepackage{latexsym}
\usepackage{amsmath}
\usepackage{amssymb}

\begin{document}
\setcounter{page}{249}
\issue{LXV~(2005)}

\title{Super Rough Semantics}

\author{A. Mani\thanks{The present author would like to thank Dr. Mohua Banerjee and Dr. J. Sen for discussions on an earlier version of the main theorem and the anonymous referee for useful comments which led to substantial improvement of the present paper.}\\
Member, Calcutta Mathematical Society\\
9B, Jatin Bagchi Road\\
Kolkata(Calcutta)-700029, India\\
\texttt{$a\_mani\_sc\_gs@yahoo.co.in$}}

\maketitle

\runninghead{A. Mani}{Rough Semantics}

\begin{abstract}
In this research a new algebraic semantics of rough set theory including additional meta aspects is proposed. The semantics is based on enhancing the standard rough set theory with notions of 'relative ability of subsets of approximation spaces to approximate'. The eventual algebraic semantics is developed via many deep results in convexity in ordered structures. A new variation of rough set theory, namely 'ill-posed rough set theory' in which it may suffice to know some of the approximations of sets, is eventually introduced.
\end{abstract}


\section{Introduction}
The approximation semantics in rough set theory is essentially captured in terms of the standard 'lower' and 'upper' approximation operators in the different forms of algebraic semantics [2, 5, 9, 13, 25, 27] known. The paper [3] contains an overview of the different approaches. None of these approaches have been used to characterize the semantics from the viewpoint of the sets of mutually approximate sets of elements derived from the original approximation space or from the viewpoint of 'the ability of subsets of the approximation space to approximate'. If $A$ is an exact element in the power set of the approximation space then it is not merely an exact element with respect to $A$. So if an element is exact then it will possess an ability to approximate. It is clear from the literature that this may possibly be expressible in terms of the topology of the rough set algebra. However this aspect does not seem to have been considered with serious practical consequences. The required 'higher order' approach is also strongly hindered by technical difficulties. In this research paper, a new semantics based on enhancing rough set theory with such a notion is developed over the concept of a 'rough algebra' originating in [2].

In the approach some meta-theoretical aspects (with respect to the approach in [2], for example) are 'internalized' and other rough theoretical concepts are supplemented. A 'definable subset' for example can be seen as a set which 'approximates' in a unique manner like no other subset. These aspects will be considered further in the last section. Some of the essential notions and terminology are repeated below for convenience. 

An \emph{approximation space} $X=\langle\underline{X},R\rangle$ is a pair with $\underline{X}$ being a set in ZFC and $R$ an equivalence relation. $\underline{X}$ can for example be taken to be a set of 'objects' and $R$ a relation which essentially assigns 'mutually exclusive types' to these objects. So $X$ is a 'qualitative structure' in the sense that the set of $R$-related objects form a lattice. For any $A\in\wp(X)$, the \emph{lower approximation} $A^l$ of $A$ is defined via,\[A^l =\bigcup\{Y; Y \in {X\backslash R},Y\subset A\}\] ($X\backslash R$ being the set of classes of $R$), while its \emph{upper approximation} $ A^u$ is defined via,\[A^u =\bigcup\{Y; Y\in X\backslash R, Y\cap A\neq\phi\}.\] $A^l$ (resp. $A^u$) can be seen to be the collection of objects that \emph{R-definitely} (resp. \emph{R-possibly}) belong to $A$ or as the collection of objects of $X$ whose types are fully included in (resp. intersect), the set of types of objects of $A$. The triple $(X,R,A)$ is called a \emph{rough set}. $A$ is \emph{roughly included} in $B$, $ A,B \subset X$, $ A \prec B$ iff $A^l\subset B^l$ and $A^u\subset B^u$. $A$ and $B$ are \emph{roughly equal} iff $A\prec B$ and $B\prec A$ iff $A^l=B^l$ and $A^u = B^u$.

A \emph{pre-rough algebra} is an algebra of the form $S\,=\,\left\langle\underline{S}, \sqcap, \sqcup, \Rightarrow, L, \neg, 0, 1, (2,2,2,1,1,0,0) \right\rangle$ which satisfies
\begin{itemize}
\item{$\left\langle\underline{S},\sqcap,\sqcup,\neg\right\rangle$ is a de Morgan lattice.}
\item{$\neg\neg{a}\,=\,a\ ;\, L(a)\sqcap{a}\,=\,L(a)$}
\item{$LL(a)\,=\,L(a)\ ;\ L(1)\,=\,1\ ;\ L(a\sqcap{b})\,=\,L(a)\sqcap L(b)$}
\item{$\neg{L}\neg{L}(a)\,=\,L(a)\ ;\ \neg{L}(a)\sqcup{L}(a)\,=\,1$}
\item{$L(a\sqcup{b})\,=\,L(a)\sqcup{L(b)}$}
\item{$(L(a)\sqcap{L(b)}\,=\,L(a),\neg{L(\neg(a\sqcap{b}))}\,=\,\neg{L}(\neg{a})\longrightarrow{a\sqcap{b}\,=\,a})$}
\item{${a\Rightarrow{b}}\,=\,(\neg{L}(a)\sqcup{L}(b))\sqcap(L(\neg{a})\sqcup\neg{L}(\neg{b}))$}
\end{itemize}
A completely distributive pre-rough algebra is called a \emph{rough algebra}. In all these algebras it is possible to define an operation $M$ by setting $M(x)\,=\,\neg{L}\neg(x)$ for each element $x$. $M$ corresponds to the upper approximation operator. The operation $\Rightarrow$ is a weaker than classical implication and corresponds to rough inclusion by way of $a\Rightarrow{b}\,=\,1$ iff $a\,\leq\,b$ in the associated lattice order. Bi-implication naturally corresponds to rough equality.

Let $H\,=\,\left\langle \underline{H},\wedge,\vee\right\rangle$ be a lattice and $T$ a binary reflexive and symmetric relation on it which is 'compatible' in the sense \[\left((a,b),(c,e)\in{T}\longrightarrow{(a\wedge{c},b\wedge{e}),(a\vee{c},b\vee{e})\in{T}}\right)\] then $T$ is called a \emph{tolerance} on $H$. A subset $B\subseteq{H}$ is called a \emph{block} of $T$ if it is a maximal subset satisfying $B^{2}\subseteq{T}$. For any $x\in{H}$ the T-associates of $x$ is the set $[x]_{T}\,=\,\{a: (a,x)\in{T}\}$. These notions directly extend to other universal algebras. A sublattice $Z$ of $H$ is called a \emph{convex sublattice} if and only if it satisfies $(\forall{x,y}\in{Z})\,x\leq{a}\leq{y}\longrightarrow {a}\in{Z}$. If $C$ is a subset of $H$ then $\downarrow{C},\uparrow{C}$ will respectively denote the lattice-ideal and filter generated by $C$. Tolerances can be fully characterized by their associated system of all blocks [7], [8] and this is denoted by $H\backslash{T}$. For finite lattices the result improves to the one presented in [7],
\begin{theorem}
If $H$ is a finite lattice, then a collection $\mathcal{S}\,=\,\{B_{\alpha}:\alpha\in{I} \}$ of subsets of $H$  is such that $\mathcal{S}\,=\,H\backslash{T}$ if and only if
\begin{enumerate}
\item{Every element of $\mathcal{S}$ is a convex sublattice of $H$.}
\item{$\mathcal{S}$ covers $H$.}
\item{$(\forall{C,E}\in{\mathcal{S}})\left(\downarrow{C}=\downarrow{E}\Longleftrightarrow{\uparrow{C}=\uparrow{E}}\right).$}
\item{For any two elements ${C,A}\in{\mathcal{S}}$ there exist ${E,F}$  such that $\left(\downarrow{C}\vee\downarrow{A}\right)\,=\,\downarrow{E},$\\ $(\uparrow{C}\vee\uparrow{A})\leq\uparrow{E},$ $\downarrow{F}\,\leq\,(\downarrow{A}\wedge\downarrow{C}),$ and $(\uparrow{C}\wedge\uparrow{A})\,=\,\uparrow{F}).$}
\end{enumerate}
\end{theorem}

A lattice is said to be \emph{semi-join distributive } if it also satisfies $(x\vee{y}\,=\,x\vee{z}\longrightarrow\,x\vee{y}\,=\,x\vee(y\wedge{z}))$. $J(L)$ will denote the set of all join-irreducible elements of a lattice $L$. A lattice $L$ is said to be \emph{finitely spatial} (resp.\emph{ spatial}) if any element of $L$ is a join-irreducible element (resp. complemented join-irreducible element) of $L$. A lattice $L$ is said to be \emph{lower continuous} if $(a\vee\,\bigwedge{X})\,=\,\bigwedge(a\vee{x})$ holds for all downward directed subsets for which $\bigwedge{X}$ exist.

In a poset $P$ a finite sequence of elements which are comparable with their predecessors is said to be a \emph{path}. A path $(x_{n})$ is said to be \emph{oriented} if $(\forall{i})\,x_{i}\leq\,x_{i+1}$ or its converse holds. A poset is said to be \emph{tree-like} if the following \emph{conditions} hold, 
\begin{itemize}
\item {If $a\leq\,{b}$ then there exists an integer $n<\,\omega$ and $x_{0},x_{1},\ldots {x_{n}}\in{P}$ such that $a\,=\,x_{0}\prec{x_{1}}\prec \ldots\prec{x_{n}}\,=\,b.$}
\item{For any two elements in the poset there exists at most one maximal path from one to the other. }
\end{itemize}
A lattice $L$ is said to be \emph{sectionally complemented} if and only if for any $b\,<\,a$ there exists a $c$ such that $b\wedge{c}\,=\,0$ and $b\vee{c}\,=\,a$.

The set of all convex sublattices of a lattice in particular and a poset in general can be endowed with a lattice structure (w.r.t inclusion) with meet corresponding to set-intersection and join corresponding to $A\vee{B}\,=\,A\,\cup\,{B}\,\cup\,\left\lbrace x:(\exists\, (y,z)\in{A\times{B}}\cup{B\times{A}}),y\,\leq\,{x}\,\leq\,{z}\right\rbrace $. The lattice is algebraic, atomistic, bi-atomic and join-semi distributive. The sublattices of such lattices have been recently characterized through [1, 29, 30, 31]. On a lattice $L$ the following conditions will be abbreviated for convenience, 
\begin{itemize}
\item {$\mathsf{S}\,:\,$  $a\wedge(b^{*}\vee{c})\,=\,(a\wedge{b^{*}})\vee{\bigvee_{i<2}}\left(a\wedge(b_{i}\vee{c})
\wedge((b^{*}\wedge(a\vee{b_{i}})\vee{c})\right)$ where $b^{*}\,=\,b\wedge(b_{0}\vee{b_{1}})$ }
\item{$\mathsf{B}\,:\,$   $ x\wedge(a_{0}\vee{a_{1}}) \wedge(b_{0}\vee{b_{1}})\,=\,\bigvee_{i\,<\,2}\left( (x\wedge{a_{i}} \wedge(b_{0}\vee{b_{1}})) \vee (x \wedge {b_{i}} \wedge (a_{0}\vee {a_{1}})\right)\vee$ \\
$ \bigvee_{i\,<\,2} \left( x \wedge (a_{0}\vee {a_{1}})\wedge (b_{0}\vee {b_{1}}) \wedge (a_{0} \vee {b_{i}})
\wedge (a_{1}\vee {b_{i-1}})\right)$}
\item{$\mathsf{U}\,:\,$   $x \wedge (x_{0}\vee {x_{1}}) \wedge (x_{1}\vee{x_{2}}) \,=\, ( x \wedge {x_{0}} \wedge (x_{1}\vee{x_{2}})) \vee \left( x \wedge {x_{1}} \wedge(x_{0}\vee{x_{2}}) \right)$\\    $ \vee \left(x \wedge{x_{2}} \wedge (x_{0}\vee{x_{1}}) \right) $}
\end{itemize}
These conditions necessarily hold in the lattice of convex sublattices of any poset. In [29, 30, 31], the following three results are proved.
\begin{theorem}
If a lattice $L$ satisfies the condition $\mathsf{S}$ then it also satisfies the condition,\\
$\mathsf{S}_{i}\,:\,(\forall{a, b, b_{0},b_{1}, c}\in{J(L)})\,$($ a\,\leq\,b\vee{c},\,b\,\leq\,b_{0}
\vee{b_{1}},\,a\neq{b}\longrightarrow $\\
$( (\exists\, \overline{b})\,a\,\leq\,\overline{b}\vee{c},\,
\overline{b}\,<\,b)\, \mathrm{or}\, ( b\,\leq\,(a\vee{b_{i}}),\,a\,\leq\,(b_{i}\vee{c})\,\mathrm{for \;some\;} \,i\,<\,2$).
\end{theorem}

\begin{theorem}
If a lattice $L$ satisfies the conditions ${\mathsf{B},\,\mathsf{U}}$ then it also satisfies the conditions $\mathsf{B}_{i}$ and $\mathsf{U}_{i}$ defined below.\\
$\mathsf{B}_{i}\,:\, (\forall{x, a_{0}, a_{1}, b_{0}, b_{1}}\in{J(L)})$($ x\,\leq\,a_{0}\vee{a_{1}},
\,b_{0}\vee{b_{1}}\longrightarrow{x}\,\leq\,a_{i},\;\mathrm{or}\;x\,\leq\,b_{i}\,\mathrm{for}\,\mathrm{some}\,i\,<\,2  $ \\
$\mathrm{or}\,x\,\leq\,a_{0}\vee{b_{0}},\,a_{1}\vee{b_{1}},\;\mathrm{or}\;x\,\leq\,a_{0}\vee{b_{1}},\,a_{1}\vee{b_{0}}$).\\
$\mathsf{U}_{i}\,:\,  (\forall{x,x_{0},x_{1},x_{2}}\in{J(L)})\left(x\,\leq\,x_{0}\vee{x_{1}},\,x_{0}\vee{x_{2}},\,
x_{1}\vee{x_{2}}\longrightarrow\,x\,\leq\,x_{0}\;\mathrm{or}\;x\,\leq\,x_{1}\;\mathrm{or}\;x\,\leq\,x_{2}
 \right) $.
\end{theorem}

If a lattice satisfies the condition $\mathsf{P}$ below, then it is said to be \emph{dually 2-distributive}.\\
$\mathsf{P}\,:\,  a\wedge(b\vee{c}\vee{e})\,=\,(a\wedge(b\vee{c}))\vee(a\wedge(b\vee{e}))\vee(a\wedge
(c\vee{e}))$.

\begin{theorem}
If $L$ is a \emph{complete, lower continuous, finitely spatial and dually 2-distributive} lattice and if it satisfies the conditions $\mathsf{S}_{i},\,\mathsf{B}_{i},\,\mathrm{and}\,\mathsf{U}_{i} $, then it satisfies the conditions $\mathsf{S},\,\mathsf{B},\,\mathrm{and}\,\mathsf{U}$.
\end{theorem}

For convenience we will call a lattice of the above form which satisfies the three conditions a \emph{long lattice}.

For more on the structure of the lattice of convex sublattices and the lattice of intervals of a lattice the reader is referred to [1, 4, 29, 30, 31]. It is proved in [16, 17] that two lattices have isomorphic convex sublattices if and only if they have isomorphic interval lattices. The following result for posets naturally applies to lattices and is relevant for the main duality result.

A poset $\Phi $ is said to be \emph{convexly isomorphic} (resp. \emph{interval isomorphic}) to another poset $\Lambda$ iff $\mathcal{F}_{cnv}(\Phi)\simeq{\mathcal{F}_{cnv}(\Lambda)}$
(resp.$\,\mathrm{ Int}(\Phi)\simeq\mathrm{Int}(\Lambda)$). $\mathcal{F}_{cnv}(\Lambda)$ being the set of all convex sublattices and $\mathrm{Int}(\phi)$ that of convex intervals of $\Lambda$ respectively. By a recent result due to [16, 17], every such $\Lambda$ is constructible from $\Phi$. One of the main results therein is stated below.

\begin{theorem}
Let $A = (\underline{A},\leq)$ be any poset. Posets \emph{convexly isomorphic} to $A$ are (up to isomorphism) just those which can be obtained by applying the following three constructions successively,
\begin{enumerate}
\item{We construct $A_1 = (\underline{A},\leq_{1})$, where $x\leq_{1}y$ means $x < y\, \mathrm{and}\,{(x,y)\notin{P}}\subset\{(x,y);x,y\in{A},x\prec{y},x\in{\mathsf{Min}(A)},y\in{\mathsf{Max}(A)}\}$.}
\item{Having $A_1$, we construct $A_2 = (\underline{A},\leq_{2})$; where $x\leq_{2}y$ holds whenever $x,y\in{C},\,x\leq_{1}y$ or $x,y\in{D},\,x\leq_{1}y$ holds for a decomposition  $A = C\cup D$ of $A$ under $(\forall{c}\in{C},d\in{D})\,c\parallel_{1}d.$ Here $\parallel_{1}$ indicates the non comparability of the two elements with respect to the order $\leq_{1}$.}
\item{Taking $A_2$, we construct $A_3=(\underline{A},\leq_{3})$; where $x\leq_{3}y$ if and only if  $x\leq{2}y$, or $(x,y)\in{Q}$ for a $Q\subset\{(x,y);(x,y)\in{A^2},\,x{\parallel_2} y,\,x\in \mathsf{Min}(A_2),\,y\in \mathsf{Max}(A_2),\}$, such that $\,(u,v),(v,w)\in Q$ do not hold simultaneously for any  $u,v,w\in{A}$.}
\end{enumerate}
\end{theorem}

\section{Main Theorems}

The new \emph{super rough algebra} will be formulated directly after a few initial definitions. The main representation theorem is proved next.
\begin{definition}
In a rough algebra of the form $S\,=\,\left\langle\underline{S},\sqcap,\sqcup,\Rightarrow,L,\neg,0,1,(2,2,2,1,1,0,0)\right\rangle$, let  $T$ be a binary relation defined by $(a,b)\in{T}$ if and only if $(\exists\, {c}\in{S}) \,L(c)\leq{a}\leq{M(c)},\,L(c)\leq{b}\leq{M(c)}$. $T$ will be called the \emph{coapproximability relation} on $S$.
\end{definition}

\begin{proposition}
The \emph{coapproximability relation} on the rough algebra is a compatible tolerance.
\end{proposition}
\begin{proof}
Clearly $T$ is both reflexive and symmetric. Let $(a,b),(c,e)\in{T}$ then there exists ${x,\,y}\in{S}$ such that the defining condition is satisfied. It follows that ${L(x\sqcap{y})}\,\leq\,{a\sqcap{c}}\leq{M(x\sqcap{y})}$ and similarly for $b\sqcap{e}$. So $(a\sqcap{c},b\sqcap{e})\in{T}$. Similarly the compatibility for $\sqcup$ holds. As the operator $L$ preserves both $\sqcap$ and $\sqcup$, so $(a,b)\in{T}\longrightarrow \left(L(a),L(b)\right)\in{T}$.

For compatibility with $\neg$, if $(a,b)\in{T}$ then there exists ${c}\in{S}\,L(c)\leq{a}\leq{\neg{L}(\neg{c})}$ and similarly for $b$. This implies $L(\neg{c})\leq\neg{a}\leq\neg{L}(c)$ and so $L(\neg{c})\leq\neg{a}\leq\neg{L}\neg(\neg{c})$ follows.
\end{proof}

With the above tolerance $T$, we can associate the collection $\mathcal{S}$ of blocks derived from the rough algebra by the basic representation theorem for tolerances on lattices. If we assume that the rough algebra is finite then the last two conditions in the representation theorem can be relaxed. From this point the choice of the 'best' underlying set for the desired representation becomes dependent on one's model-theoretical preferences. The first option is to take the collection of all convex sublattices of the rough algebra as the underlying set, define other 'global operations', and specify a scheme for deriving the system of blocks. The second option is to start from the set of blocks, adjoin the set of principal ideals and filters generated by it, 'partially complete' their 'join' and 'meet' and then define the required global operations. A third option is to take the union of the set of blocks, the set of all filters and ideals as the underlying set. All of these result in different amounts of partiality of some of the operations. In the first approach the only partial operation is the operation $L_{T}$, the partiality is not particularly significant.  For more on partial algebras [6, 18] can be consulted for example. \\

\begin{proposition}
Every block of the tolerance $T$ is an interval of the form $\left[ \bigwedge{a_{i}},\bigvee{a_{i}}\right] $, whenever the approximation space is finite.
\end{proposition}
\begin{proof}
If $B$ is a block of a tolerance of a lattice then it's supremum and infimum must necessarily exist by the result in [7]. So the result follows with the elements $a_{i}$ being the elements of the block.
\end{proof}

In a partial algebra of the form $S\,=\,\left\langle\underline{S},f_{1},f_{2},\ldots\right\rangle $, two term operations $t,\,s$ are said to be weakly equal in $S$, i.e. $\left( t^{\underline{S}}\,\stackrel{\omega}{=}\,s^{\underline{S}}\right) $ if and only if, $(\forall{x}\in\,{\mathsf{dom}(t^{\underline{S}})\,\cap\,\mathsf{dom}(s^{\underline{S}})})\,\left( t^{\underline{S}}(x)\,{=}\,s^{\underline{S}}(x)\right) $ (i.e. if the left and the right hand side are defined, then the two are equal). Actually weak equations are special types of ECE-equations (see [6]). The strong weak equality $\stackrel{\omega^{*}}{=}$ is defined in the same way with the quantification part being modified to $(\forall{x}\in\,{\mathsf{dom}(t^{\underline{S}})\,=\,\mathsf{dom}(s^{\underline{S}})})$, (i.e. if either side is defined then the other is, and the two are equal). Since the interpretation part used in the partial algebras will be quite direct, so we will use relatively simplified notation.

In the following, except for $\sqcap,\,\sqcup\,\mathrm{and}\,L$ we will use the same operation symbol to denote the operation itself in the rough algebra and the super rough set-algebra. The interpretation should be clear from the context. For example, $\sqcap^{\underline{\Re}}$ means the interpretation of the operation symbol $\sqcap$ over $\underline{\Re}$.

\begin{definition}
A \emph{super rough set-algebra} will be a partial algebra of the form       \[ \Re\,=\,\left\langle\underline{\Re},\wedge,\vee,\sqcap^{\underline{\Re}},\sqcup^{\underline{\Re}},\neg,L^{\underline{\Re}},L_{T}, \downarrow, \uparrow, \underline{S},\emptyset,(2,2,2,2,1,1,1,0,0)\right\rangle\]
satisfying all of the following conditions:
\begin{enumerate}
\item{The underlying set of the finite rough algebra $S$ (with the above indicated operations) is $\underline{S}$.}
\item {The set of all convex sublattices of $S$ is $\underline{\Re}$.}
\item{$(\forall{A,B}\in{\Re})\,A\sqcap^{\underline{\Re}}{B}\,=\,\left\{x\sqcap{y}:x\in{A},\,y\in{B}\right\}$ if defined in $\Re$.}
\item{$(\forall{A,B}\in{\Re})\,A\sqcup^{\underline{\Re}}{B}\,=\,\left\{x\sqcup{y}:x\in{A},\,y\in{B}\right\}$ if defined in $\Re$.}
\item{The principal filter and the principal ideal operations with respect to the lattice order on the set of convex sublattices are respectively $\uparrow$ and $\downarrow$.}
\item{The usual lattice operations on the set of all convex sublattices will correspond to $\wedge$ and $\vee$.}
\item{$\neg{A}\,=\,\left\{\neg{x}:\,x\in{A}\right\}.$}
\item{$L^{\underline{\Re}}(A)\,=\,\left\{L(x):x\in{A}\right\}$ if defined in $\Re$.}
\item{$\Re\,\models\,\mathsf{S}_{i},\,\mathsf{B}_{i},\,\mathsf{U}_{i}$}
\item{\flushleft{\[L_{T}(A)\,=\,\left\{\begin{array}{ll}      A & {\mathrm{if} \,A\,\mathrm{  is\,  a\,  block\,  of\,}  {T}.} \\
 {\mathrm{undefined}} & {\mathrm{else}}
\end{array}\right.\]}}\end{enumerate}
\end{definition}
\begin{theorem}
A super rough set-algebra $\Re$ as defined above satisfies all of the following: \begin{enumerate}
\item{$\left\langle \underline{\Re},\sqcup^{\underline{\Re}},\sqcap^{\underline{\Re}}\right\rangle$ is a  partial distributive lattice which satisfies all weak-equalities in the $\omega^{*}$ sense.}
\item{$\wedge,\vee$ are total lattice operations on $\Re$.}
\item{$\neg(A\sqcap^{\underline{\Re}}{B})\,=\,(\neg{A}\sqcup^{\underline{\Re}}\neg{B}).$}
\item{$\neg\neg{A}\,=\,{A}$}
\item{$L^{\underline{\Re}}(A)\,\stackrel{\omega^{*}}{=}L^{\underline{\Re}}(L^{\underline{\Re}}(A))$}
\item{$L^{\underline{\Re}}(A)\,\stackrel{\omega^{*}}{=}L^{\underline{\Re}}(A)\sqcap{A}$}
\item{$L^{\underline{\Re}}(A\sqcap{B})\,\stackrel{\omega^{*}}{=}L^{\underline{\Re}}(A)\sqcap{L^{\underline{\Re}}(B)}$}
\item{$L^{\underline{\Re}}(A\sqcup{B})\,\stackrel{\omega^{*}}{=}L^{\underline{\Re}}(A)\sqcup{L^{\underline{\Re}}(B)}$}
\item {$(L_{T}(A)\,=\,A \longrightarrow\,(\{x\}\wedge{A}\,=\,\{x\}\longleftrightarrow \{x\}\sqcup\{\bigwedge{A}\},\,=\,\{x\},\,\{x\}\sqcap\{\bigvee{A}\}\,=\,\{x\})) $}
\item {$(\forall{A,\,B})(L_{T}(A)\,=\,A, L_{T}(B)\,=\,B,\, (\uparrow{A})\,=\,(\uparrow{B}) \longrightarrow\,{(\downarrow{A})\,=\,(\downarrow{B})})$}
\item {$(\forall{A,\,B})(L_{T}(A)\,=\,A, L_{T}(B)\,=\,B,\,(\downarrow{A})\,=\,(\downarrow{B}) \longrightarrow\,{(\uparrow{A})\,=\,(\uparrow{B})})$}
\item{$(L_{T}(A)\,=\,A\longrightarrow\,L(A)\vee{A}\,=\,A,\neg{L}(\neg{A})\vee{A}\,=\,A )$}
\item{$(L_{T}(A)\,=\,A\longrightarrow\,L(A)\sqcup{A}\,=\,A,\neg{L}(\neg{A})\sqcup{A}\,=\,A)$ }
\item {For any two fixed points $A,\,B$ of $L_{T}$, there exist two other fixed points $E,\,F$ such that\\ $((\downarrow{A})\vee(\downarrow{B}))\,=\,(\downarrow{E})$, $((\uparrow{B})\vee(\uparrow{A}))\leq(\uparrow{E})$, $(\downarrow{F})\leq((\downarrow{A})\wedge(\downarrow{B}))$ and  $((\uparrow{B})\wedge(\uparrow{A}))\,=\,(\uparrow{F})$ hold.}
\item {$(L_{T}(A)\,=\,A\longrightarrow\,(\exists {B})\,\neg(A)\wedge{B}\,=\,\neg{A},\,L_{T}(B)\,=\,B)$}
\item {$(L_{T}(A)\,=\,A\longrightarrow\,(\exists {B})\,L(A)\wedge{B}\,=\,L(A),\,L_{T}(B)\,=\,B)$}
\item {$(L_{T}(A)\,=\,A,\,L_{T}(B)\,=\,B,\,A\sqcap{B}\,=\,C\longrightarrow\,(\exists{E})\,E\wedge{C}\,=\,C, \,L_{T}(E)\,=\,E)$}
\item {$(L_{T}(A)\,=\,A,\,L_{T}(B)\,=\,B,\,A\vee{B}\,=\,C\longrightarrow\,(\exists {E})\,E\wedge{C}\,=\,C, \,L_{T}(E)\,=\,E)$}
\item {$(L_{T}(A)\,=\,A,\,L_{T}(B)\,=\,B,\,A\wedge{B}\,=\,C\longrightarrow\,(\exists {E})\,E\wedge{C}\,=\,C, \,L_{T}(E)\,=\,E)$}
\item {$(L_{T}(A)\,=\,A,\,L_{T}(B)\,=\,B,\,A\sqcup{B}\,=\,C\longrightarrow\,(\exists {E})\,E\wedge{C}\,=\,C, \,L_{T}(E)\,=\,E )$}
 The order relation used is the one on the lattice of convex sublattices.
\end{enumerate}
\end{theorem}
\begin{proof}
\begin{enumerate}
\item {To prove $\left\langle \underline{\Re},\sqcup,\sqcap\right\rangle$ is a  partial distributive lattice which satisfies all weak-equalities in the $\omega^{*}$ sense, we need to prove the strong weak associativity, commutativity, idempotence, distributivity and then prove absorption. Idempotence obviously holds in the total sense as $A\sqcup{A}$ is always defined and equal to $A$. For distributivity, if $x\in{A\sqcap(B\sqcup{C})}$ then there necessarily exist $a\in{A},\,b\in{B},\,c\in{C}$ such that $x\,=\,a\sqcap(b\sqcup{c})\,=\,(a\sqcap{b})\sqcup(a\sqcap{c})$. So $x\in {(A\sqcap{B})\sqcup(A\sqcap{B}}$ and vice versa. Absorption is true in the sense $(A\sqcap({B}\sqcup{A})\,=\,C\longrightarrow {A\,=\,C})$ and it's dual as we are dealing with blocks. }
\item {It is known that $\wedge,\vee$ are total lattice operations on $\Re$.}
\item {If $x\in{\neg(A\sqcap{B})}$, then there exist $a\in{A},\,b\in{B}$ such that $x\,=\,\neg(a\sqcap{b})\,=\,(\neg{a}\sqcup\neg{b})$. So $x\in{(\neg{A}\sqcup\neg{B})}$ and conversely. }
\item {For any $A$, the equality $\neg\neg{A}\,=\,{A}$ is obvious.}
\item {If $L(A)$ is defined, then it is necessary that $L(L(A))$ be defined and the two must be equal. The converse also holds.}
\item {If $L(A)$ is defined then it is necessarily a convex subset and $L(A)\sqcap{A}\,=\,\left\lbrace L(x)\sqcap{y}\,:\,x,\,y\in{A}\right\rbrace$. Clearly $\left( L(x)\sqcap{L(y)}\,\leq\,L(x)\sqcap{y}\,\leq\,L(x),y\right) $ in the convex set (the order $\leq$ being the one on the underlying rough algebra). So $L(x)\sqcap{y}$ must be in $L(A)$. }
\item{A two-way inclusion argument suffices for the condition,  $L^{\underline{\Re}}(A\sqcap{B})\,\stackrel{\omega^{*}}{=}L^{\underline{\Re}}(A)\sqcap{L^{ \underline{\Re}}(B)}$. The definability part is direct.}
\item {The proof of $L^{\underline{\Re}}(A\sqcup{B})\,\stackrel{\omega^{*}}{=}L^{\underline{\Re}}(A)\sqcup{L^{ \underline{\Re}}(B)}$ is as in the above.}
\item {The proof of $(L_{T}(A)\,=\,A \longrightarrow\,(\{x\}\wedge{A}\,=\,\{x\}\longleftrightarrow \{x\}\sqcup\{\bigwedge{A}\},\,=\,\{x\},\,\{x\}\sqcap\{\bigvee{A}\}\,=\,\{x\})) $ follows from Proposition 2.2 above. }
\item {The representation theorem for tolerances in finite lattices includes $(\forall{A,\,B})(L_{T}(A)\,=\,A, L_{T}(B)\,=\,B,\,\uparrow{A}\,=\,\uparrow{B}\longrightarrow\, {\downarrow{A}\,=\,\downarrow{B}})$.}
\item {The proof of $(\forall{A,\,B})({L_{T}(A)\,=\,A, L_{T}(B)\,=\,B},\downarrow{A}\,=\,\downarrow{B} \longrightarrow\,\uparrow{A}\,=\,\uparrow{B})$ is as in the above.}
\item {If $L_{T}(A)\,=\,A$, then $A$ is necessarily a block of the tolerance $T$. For $(L_{T}(A)\,=\,A\longrightarrow\,L(A)\vee{A}\,=\,A,\neg{L}(\neg{A})\vee{A}\,=\,A )$. If $A$ is a block, then $L(A)$ is also convex. In the convex order, $A\vee{L(A)}\,=\,A\,\cup\,{L(A)}\,\cup\, \left\lbrace x:\exists\, (y,z)\in{A\times{L(A)}}\cup{L(A)\times{A}},y\,\leq\,{x}\,\leq\,{z}\right\rbrace $. But $x,y\in{A}$ implies $(x,y)\in{T}$. If $(x,y)\in{T}$ and $\prec$ is the lattice order on the rough algebra, then ${L(x)}\prec{x}$ and $(\exists a),\,L(a)\prec{x},y\prec{\neg{L}\neg{a}}$. Suppose $L(a)\parallel{L(x)}$, then $L(a\sqcap(x\sqcup{y}))\prec{L(x)\prec\,{\neg{L}\neg(a\sqcap(x\sqcup{y}))}}.$ Hence $(L(x),y)\in{T}$. This ensures the conclusions. }
\item {$(L_{T}(A)\,=\,A\longrightarrow\,L(A)\sqcup{A}\,=\,A,\neg{L}(\neg{A})\sqcup{A}\,=\,A)$ is direct. }
\item {This is essentially the last condition in the representation theorem by blocks for tolerances of a finite lattice.  }
\item {For  $(L_{T}(A)\,=\,A\longrightarrow\,(\exists {B})\,\neg(A)\wedge{B}\,=\,\neg{A},\,L_{T}(B)\,=\,B)$, note that if $A$ is a fixed point of $L_{T}$ then it is a block. The set $\neg{A}$ will then be such that if $\alpha\,\in\,(\neg{A})^{2}$ then $\alpha\,\in\,T$ by the compatibility of the associated operation with the tolerance. But by the representation theorem for tolerances on an algebra [7], there must exist a block $B$ containing $\neg{A}$.}
\end{enumerate}
The rest follow by considerations similar to that used for the last item.
\end{proof}

Based on the above concrete situation we introduce the concept of a \emph{Super Rough Algebra}.

\begin{definition}
A partial algebra of the form $S\,=\,\left\langle \underline{S},\wedge, \vee, \sqcap, \sqcup, \neg, L, L_{T}, \downarrow, \uparrow, 1, \emptyset \right\rangle$ with associated arities $(2,2,2,2,1,1,1,0,0)$, will be called a \emph{super rough algebra} if all of the following hold : \begin{enumerate}
\item {$S\,=\,\left\langle \underline{S},\wedge, \vee,\right \rangle$ is a long lattice.}
\item{$\neg(a\sqcap{b})\,=\,(\neg{a}\sqcup\neg{b}).$}
\item{$\neg\neg{a}\,=\,{a}$}
\item{$L(a)\,\stackrel{\omega^{*}}{=}L(L(a))$}
\item{$L(a)\,\stackrel{\omega^{*}}{=}L(a)\sqcap{a}$}
\item{$L(a\sqcap{b})\,\stackrel{\omega^{*}}{=}L(a)\sqcap{L(b)}$}
\item{$L(a\sqcup{b})\,\stackrel{\omega^{*}}{=}L(a)\sqcup{L(b)}$}
\item{$(L_{T}(a)\,=\,a\longrightarrow\,L(a)\vee{a}\,=\,a,\neg{L}(\neg{a})\vee{a}\,=\,a )$}
\item{$(L_{T}(a)\,=\,a\longrightarrow\,(\{x\}\wedge{a}\,=\,\{x\}\longleftrightarrow\,\{x\}\sqcup\{
\bigwedge{a}\}\,=\,\{x\},\{x\}\sqcap\{\bigvee{a}\}\,=\,\{x\}))$}
\item{$(L_{T}(a)\,=\,a\longrightarrow\,L(a)\sqcup{a}\,=\,a,\neg{L}(\neg{a})\sqcup{a}\,=\,a )$}
\item{$(\forall{x,a})\,( a\wedge{x}\,=\,a,\,a\wedge{\emptyset}\,=\,\emptyset \longrightarrow\,a\,=\,x,(\exists\, {y})\,L_{T}(y)\,=\,y,\,y\wedge{x}\,=\,x)$}
\item{$(\forall{x,y})\,\left(L_{T}(x)\,=\,x,\,L_{T}(y)\,=\,y,(\downarrow{x})\,=\,(\downarrow{y})\longrightarrow{(\uparrow{x})\,=\,(\uparrow{y})}\right)$}
\item{$(\forall{x,y})\,\left(L_{T}(x)\,=\,x,\,L_{T}(y)\,=\,y,(\uparrow{x})\,=\,(\uparrow{y})\longrightarrow{(\downarrow{x})\,=\,(\downarrow{y})}\right)$}
\item{For any two fixed points $a,\,b$ of $L_{T}$, there exist two other fixed points $x,\,y$ such that \\ $((\downarrow{a})\vee(\downarrow{b}))\,=\,(\downarrow{x})$, $((\uparrow{b})\vee(\uparrow{a}))\leq(\uparrow{x})$, $(\downarrow{y})\leq((\downarrow{a})\wedge(\downarrow{b}))$ and  $((\uparrow{b})\wedge(\uparrow{a}))\,=\,(\uparrow{y})$ hold.}
\end{enumerate}
\end{definition}

\begin{theorem}
For every super rough algebra $S$, there exists an approximation space $X$ such that the super rough set algebra generated by $X$ is isomorphic to $S$.
\end{theorem}
\begin{proof}
As $S$ is a long lattice so there exists a partially ordered set $P$ such that the lattice of convex subsets $Co(P)$ generated by it is isomorphic to a sublattice of it. The finiteness part ensures that we can obtain it as an isomorphic copy of $Co(P)$. Actually this part is not essential for our proof. The convex structure simply ensures better expression in terms of 'total' operations as opposed to 'partial' ones and is always available.

By the fundamental characterization of tolerances by blocks, we can immediately reconstruct
a lattice $F$ alongwith a compatible tolerance $T$ on it from the set of fixed points of the map $L_{T}$. $L_{T}$ is actually a definable via the conditional equations. Note that $F$ is also constructible as the set of 'singletons' in $S$. These are definable via the covering property with respect to the empty set. Again note that in any partially ordered set all singletons are convex subsets. This allows the definition of the operations $\sqcup,\,\sqcap,\,L,\,\neg$ and the distinguished elements on the desired prerough algebra.

Now the concrete representation theorem for prerough algebras proved in [2] allows the existence of the approximation space $X$. Checking that the super rough set algebra generated by the approximation space is isomorphic to $S$ is by a direct contradiction argument.
 
\end{proof}

\begin{remark}
The proof is substantially simplified by the two representation theorems. As the algebras are finite there is no distinction between rough and prerough algebras.
\end{remark}

\begin{remark}
The conditions defining a long lattice are not difficult to check in actual usage. See [29, 30, 31] for the details.
\end{remark}

\section{Applicability of the Coapproximability Semantics-Ill-Defined Approximation Spaces}

A natural interpretation of the above can be as an alternative algebraic semantics of rough set-theoretical reasoning. This will be called the \emph{direct alternative semantics}. But as it is formulated with many other meta-theoretical aspects included, this view is incomplete. The co-approximation can also be viewed as a density relation. But as it is a tolerance and not an equivalence this co-approximation viewpoint requires more investigation of the 'uniformities' which sustain the view.

Again since the coapproximation relation is a tolerance we can view the set of blocks as a reconstruction of the 'powerset' (derived from the approximation space)  but from the
approximation semantics perspective. So what does this mean from an application oriented viewpoint?  Simply that the semantics motivates a procedure for fitting the model to situations where we see 'approximations working'. This is vaguely similar to statistical procedures in which we assume 'randomness' or a 'specific distribution'. Here we assume there is a 'rough approximation operating on an unknown approximation space'.  By an 'unknown approximation space' of course we mean a relatively ill-defined approximation space in which not all of the approximations are reliable or known. 

Suppose we have a composite property $\pi$ for which we are to develop a rough semantics on an approximation space. Suppose further that the connection of the property $\pi$ with the type of attributes is not wholly clear (possibly in a mathematically measurable sense). For some objects this can be  acceptable. The above semantics essentially provides a means of improving on the available information using the objects for which 'it is acceptable'. Suitable \emph{acceptor functions} will be usually definable in the situation.

In the direct alternative semantics-perspective note that given a pair of sets we can decide on the existence of a lower and upper approximation pair (which includes them) and formulate a procedure for computing such definite sets.

\subsection{Extended Example}
In the following we apply the developed theory in an evaluation contest. As such it constitutes a new way of testing and evaluation in the context.

Suppose we need to test a group of students in accessing the quality of fruits of a species. We expect the students to have a good conception of the abstractable qualities of the fruits. If a student is given a set of fruits for categorization, then he/she is liable to require less conception of the abstract, than when comparing just two of the fruits. This is because there is more scope for confirming the 'interpretation in models' (or simply in application) of the abstract concepts involved.

So our test procedure will consist in giving a pair of sets of fruits to the students at a time and asking them to select just one of the following four responses :

\begin{itemize}
\item {Fruit set $A$ has all the positive features that fruit set $B$ has. }
\item {Fruit set $B$ has all the positive features that fruit set $A$ has. }
\item {Both sets of fruits $A$ and $B$ have noncommon features.}
\item {Fruit set $A$ and $B$ have identical features.}
\end{itemize}

Let $S$ and $\wp(S)$ be the set of fruits used and its power set respectively. Let $F$ be a subset of $\wp(S)$. Each student must preferably be given a set of pairs of sets $K$ (say), such that $K$ admits of a decomposition of the form $K\,=\,\bigcup_{x\in{F}}\left\lbrace (x,y)\,:y\in{S}\right\rbrace $.

Apparently, the test procedure concerns the skill level of the students in deciding the quality of fruits. But using the developed procedure we can infer a lot more about the student's  comprehension of the abstractable qualities of the fruits. For this we must simply form the super rough algebra from the given information and see which of the rough algebra semantics it corresponds to. The semantics being formulated by using different levels of abstractable qualities in the usual rough way. In the approach if a small number of bad errors are made by the students then the result can be severe. However we can use extra procedures to relax the severity if desired.

 \subsection{ Ill-Posed Problems}

Different types of ill-posed problems can be solved using the developed theory. The problem may relate to deciding in situations where we have multiple sets of inconclusive rough information (via multiple agents for example) at one level. The problem may also involve partial information structures, in the sense that the approximations following from the original approximation space may not be wholly reliable (but a reliable fragment is identifiable). Other types of problems include all those which have associated difficulties in rough approaches at the immediately higher meta-level of discourse.

Suppose we have an approximation space $X$ and let the rough algebraic semantics for it be available. Now if the problem is to reconsider the entire problem on the basis of newly available information, attributes and/or objects, then the problem of extending the semantics or formulating a new rough semantics is solvable by the above. An important subclass of problems includes those in which the extra information is of the form $\exists\, {A}\,A^{l}\leq\,B,\,C\,\leq\,A^{u}$ for given pairs $B$ and  $C$. If we try to use additional information about 'approximations' to be fitted in, then in the  usual perspective it generally means difficult embedding problems. These have also not been considered in much detail. The study of rough equalities in [21] is suggestive of the complexity of the problem. The greatest advantage over the 'rough-algebra perspective' is that we can make effective use of the knowledge of betweenness of sets within upper and lower approximations of some other set.

\subsection{Super Rough Logic}

Interpreting the super rough algebra itself as an algebra of a logic in a modified Blok-Pigozzi
or\\ Czelakowski sense we can possibly form a new logic. The logic can be seen as a logic of coapproximable entities which includes usual rough logic completely. The most interesting feature of interest in the logic would once again be the 'ability to manage with much less information'. This however need not constrain the other features arising at 'blocks'. It may be noted that rough set theories starting from generalized approximation spaces involving tolerances [32] are quite unrelated to the present questions, though the present theory can possibly be extended to those spaces.

\section*{Conclusions}

In this research paper, a new higher order semantics of rough sets is progressed. The main advantage of the semantics is in its flexibility in possible applications. In particular it allows a framework for analyzing partially defined rough information and completing them. At another level it allows for revising a given rough semantics in the light of additional rough information. The latter two problems have been termed 'ill-posed' because of their natural character in usual rough analysis. Something more needs to be done in the computational part for associated applications. We expect substantial expansion of the bounds of rough analysis by the present approach.

\begin{small}

\end{small}

\end{document}